\documentclass[11pt]{article}
\usepackage{setspace}
\usepackage{amsmath,amssymb,amsfonts,amsthm}
\newcounter{item}[section]
\newcounter{kirshr}
\newcounter{kirsha}
\newcounter{kirshb}
\newtheorem{theorem}{Theorem}[section]

\newtheorem{corollary}[theorem]{Corollary}
\theoremstyle{definition}

\newtheorem{definition}[theorem]{Definition}

\def\(R)RA{{\bf (R)RA}}

\begin{document}

Thanks to Dr. Tarek Sayed Ahmed whose recent work was a motivation of the subject of this paper.

Let $L$ denote a non-empty countable relational language (this entails no loss of generality): $L=(R_i)_{i \in I}$ where $I$ is a non-empty countable index set and $R_i$ is an $n_i$-ary relation symbol. Denote by $X_L$ the space $$X_L=\prod_{i \in I} 2^{(\mathbb{N}^{n_i})}.$$ We view the space $X_L$ as the space of countably infinite \emph{L-structures}.

A \emph{fragment} $F$ of $L_{\omega_{1}\omega}$ is a set of formulas in $L_{\omega_{1}\omega}$ containing all atomic formulas, closed under subformulas, negation, quantifiers and finite conjunctions and disjunctions.

\begin{definition}
For $\varphi(=\varphi(\bar{v}))$ a formula of $L_{\omega_{1}\omega}$ and $\bar{s}$ a finite sequence from $\mathbb{N}$ of appropriate length (i.e, $\bar{s}\in{}^{|\Delta\varphi|}\omega$), let $$Mod(\varphi,\bar{s})=\{x\in X_L: \mathcal{U}_x\models\varphi[\bar{s}]\},$$ where $\varphi[\bar{s}]$ denotes the sentence obtained from the formula $\varphi(\bar{v})$ by substituting $\bar{s}$ for the free variables. (If $\varphi$ is a sentence, we write $Mod(\varphi)$ for $Mod(\varphi,()).$).
\end{definition}

Let $t_F$ be the topology on $X_L$ generated by $\mathcal{B}_F=\{Mod(\varphi,\bar{s}):\varphi\in F, \bar{s}\in{}^{|\Delta\varphi|}\omega\}.$ By a result of Sami (See [4]), $t_F$ is a Polish topology on $X_L.$

Let $F$ be a fragment of $L_{\omega_{1}\omega}$. We say that $x,y\in X_L$ (or \emph{their corresponding structures}) are \textit{separable in} $F$, if there is $\varphi\in F$ such that $|\varphi^x|\neq|\varphi^y|,$ where $\varphi^x=\{\bar{s}\in{}^{|\Delta\varphi|}\omega:\mathcal{U}_x\models\varphi[\bar{s}]\}$. (It is clear that if two structures are separable in some fragment, then they are non-isomorphic). Notice that, if $\varphi$ is a sentence, then for all $x$, either $\varphi^x$ is empty or else contains only the empty sequence.

For $F$ a fragment
$$E_F=\{(x,y)\in X_L\times X_L: \mbox{ For all }\varphi\in F, \mbox{ }|\varphi^x|=|\varphi^y|\}.$$

\begin{theorem}
For $F$ a countable fragment of $L_{\omega_{1}\omega}$, $E_F$ is Borel in the product topology $(X_L,t_F)\times(X_L,t_F)$.
\end{theorem}
For every $\varphi\in F$ which is not a sentence, select a bijection $\mu_\varphi:\mathbb{N}\longrightarrow{}^{|\Delta \varphi|}\omega$.
If $\varphi$ is a sentence, let $\mu_\varphi$ be the constant map from $\mathbb{N}$ to $\mathbb{N}$ that sends everything to $1$ (a value that cannot be the empty sequence). It is clear that, for a set $X\subseteq {}^{|\Delta\varphi|}\omega$,
$$X\mbox{ is infinite \textit{iff}  } (\forall n)(\exists m>n) \mu_\varphi(m)\in X.$$

\begin{align*}
 |X|=|Y|\in\omega &\Longleftrightarrow (\exists n)(\exists f,g\in Inj(n,{}^{|\Delta \varphi|}\omega))(f^*(n)=X\wedge g^*(n)=Y)\\
  &\Longrightarrow (\exists n)(\exists f,g\in Inj(n,{}^{|\Delta \varphi|}\omega))(f^*(g^{-1}(Y))=X\wedge \\ &g^*(f^{-1}(X))=Y) \\
 &\Longrightarrow (\exists n)(\exists f,g\in Inj(n,{}^{|\Delta \varphi|}\omega))(X\subseteq f^*(n)\wedge Y\subseteq g^*(n)\\&\wedge g^{-1}(Y)=f^{-1}(X))\\
 &\Longrightarrow (\exists n)(\exists f,g\in Inj(n,{}^{|\Delta\varphi|}\omega))(|X|=|f^{-1}(X)|\wedge |Y|=\\&|g^{-1}(Y)|\wedge |g^{-1}(Y)|=|f^{-1}(X)|)\\
 &\Longrightarrow |X|=|Y|\in\omega.
\end{align*}

\begin{align*}
f^*(g^{-1}(Y))=X
&\Longleftrightarrow (\forall t)[t\in X\Leftrightarrow g(f^{-1}(t))\in Y].
\end{align*}

Now we are ready to prove our main theorem.

\begin{corollary}
Let $T$ be a first order theory in a countable language. If $T$ has an uncountable set of pairwise separable (in any countable fragment of $L_{\omega_1\omega}$) countable models, then it has such a set of size $2^{\aleph_0}$ (and so has $2^{\aleph_0}$ non-isomorphic countable models).
\end{corollary}

The above corallary can have other versions. We can talk about any set of models of $T$ whose corresponding set of codes is  $G_\delta$ in $X_L$. For example, suppose we are given a certain countable family, $\{\Gamma_i:i<\omega\},$ of non-isolated $n$-types ($n\in\omega$) of $T$ (see [3]).

\end{document}